\documentclass[12pt]{article}
\usepackage{amsmath,amsfonts,amssymb}

\def\x{{X}}
\def\y{{Y}}
\def\Z{\mathbb Z}
\def\lm{{\lambda/\mu}}
\def\inv{\textnormal{inv}}

\def\qed{\nopagebreak\hfill{\rule{4pt}{7pt}}}
\def\proof{{\it{Proof.} \hskip 2pt}}

\newdimen\Squaresize \Squaresize=14pt
\newdimen\Thickness \Thickness=0.7pt

\def\Square#1{\hbox{\vrule width \Thickness
   \vbox to \Squaresize{\hrule height \Thickness\vss
      \hbox to \Squaresize{\hss#1\hss}
   \vss\hrule height\Thickness}
\unskip\vrule width \Thickness} \kern-\Thickness}

\def\Vsquare#1{\vbox{\Square{$#1$}}\kern-\Thickness}

\def\moins{\raise 1pt\hbox{{$\scriptstyle -$}}}

\newtheorem{theo}{Theorem}[section]
\newtheorem{deff}[theo]{Definition}
\newtheorem{lemm}[theo]{Lemma}

\newtheorem{coro}[theo]{Corollary}

\makeatletter \@addtoreset{figure}{section} \makeatother
\makeatletter
\long\def\@makecaption#1#2{%
   \vskip 10\p@
   \setbox\@tempboxa\hbox{{#1}\ \ #2}%
   \ifdim \wd\@tempboxa >\hsize
       {#1}\ \ #2\par
   \else
       \hbox to\hsize{\hfil\box\@tempboxa\hfil}%
   \fi}
\makeatother

\begin{document}
\begin{center}
{\Large \bf Transformations of Border Strips and \\
\vspace{2mm}Schur Function Determinants}
\end{center}

\begin{center}
{William Y. C. Chen$^{1}$, Guo-Guang Yan$^{2}$ and Arthur L. B. Yang$^{3}$}\\
\vspace{2mm}

Center for Combinatorics, LPMC\\
Nankai University, Tianjin 300071, P. R. China\\
\vspace{2mm}
 Email: $^{1}$chen@nankai.edu.cn,
$^{2}$guogyan@eyou.com, $^{3}$arthuryang@nankai.edu.cn
\end{center}

\vspace{0.3cm} \noindent{\bf Abstract.} We introduce the notion of
the cutting  strip of an outside decomposition of a skew shape,
and show that cutting strips are in one-to-one correspondence with
outside decompositions for a given skew shape. Outside
decompositions are introduced by Hamel and Goulden and are used to
give an identity for the skew Schur function  that unifies the
determinantal expressions for the skew Schur functions including
the Jacobi-Trudi determinant, its dual, the Giambelli determinant
and the rim ribbon determinant due to Lascoux and Pragacz. Using
cutting strips, one obtains a formula for the number of outside
decompositions of a given skew shape. Moreover, one can define the
basic transformations which we call the twist transformation among
cutting strips, and derive a transformation theorem for the
determinantal formula of Hamel and Goulden. The special case of
the transformation theorem for the Giambelli identity and the rim
ribbon identity was obtained by Lascoux and Pragacz. Our
transformation theorem also applies to the supersymmetric skew
Schur function.

 \vskip 10pt


\noindent {\bf Suggested Running title:} Border Strips


\noindent {\bf Corresponding Author:} William Y. C. Chen,
chen@nankai.edu.cn

\section{Introduction}

We assume that the reader is familiar with the background of the
determinantal expressions for the skew Schur function
$s_{\lm}(\x)$ in the variable set $\x=\{x_1,x_2,\ldots\}$,
including the Jacobi-Trudi determinant, its dual, the Giambelli
determinant and the Lascoux-Pragacz rim ribbon determinant. Here
are a few relevant references \cite{ER, FK, GV1, GV2, LP2, Mac1,
Sagan, S1, Stem, Ueno}.

This paper is motivated by the  remarkable discovery of Hamel and
Goulden \cite{HG} that the four determinantal formulas for the
skew Schur functions can be unified through the concept of outside
decompositions of a skew shape. For any outside decomposition,
Hamel and Goulden derive a determinantal formula with strip Schur
functions as entries. Their proof is based on a lattice path
construction and the Gessel-Viennot methodology \cite{GV2}. The
first result of this paper is the notion of cutting strips of
outside decompositions for a given skew shape. We show that each
outside decomposition can be recovered from a strip, which we call
the cutting strip. It turns out that the number of outside
decompositions follows from the enumeration of cutting strips.
Moreover, the cutting strips can be used to give an easier
construction of the strips involved in the Hamel-Goulden
determinant.

The second result of this paper is to give a transformation
theorem for the Hamel-Goulden determinant. Using cutting strips to
represent outside decompositions, one can define the basic steps
to transform any outside decomposition to another which we call
the twist transformations. For the twist transformations on
outside decompositions,  we give a construction of determinantal
operations that transform the determinantal formula for one
outside decomposition to the determinantal formula corresponding
to any other outside decomposition.

We remark that the arguments in this paper can be carried over to
the case of supersymmetric skew Schur functions and the
corresponding determinantal formula.

\section{Cutting Strips}

We begin this section with a review of relevant terminology and
notation. Let $\lambda$ be a \emph{partition} of $n$ with at most
$k$ parts, i.e., $\lambda=(\lambda_1,\lambda_2,\ldots,\lambda_k)$
where $\lambda_1\geq\lambda_2\geq\ldots\geq\lambda_k\geq0$ and
$\lambda_1+\lambda_2+\ldots+\lambda_k=n$. The \emph{(Young} or
\emph{Ferrers) diagram} of $\lambda$ may be defined as the set of
points $(i,j)\in \Z^2$ such that $1\leq j\leq\lambda_i$ and $1\leq
i\leq k$. It can be represented in the plane by an array of square
boxes or cells that is top and left justified with $k$ rows and
$\lambda_i$ boxes in row $i$. A box $(i,j)$ in the diagram is the
box in row i from the top and column $j$ from the left. Given two
partitions $\lambda$, $\mu$ with at most $k$ parts, and
$\lambda_i\geq\mu_i$, for $i=1,2,\ldots,k$, then $\lm$ is a skew
partition, and the \emph{diagram of skew shape} $\lm$ or the
\emph{skew diagram} of $\lm$ is defined as the set of points
$(i,j)\in \Z^2$ such that $\mu_i< j\leq\lambda_i$ and $1\leq i\leq
k$. In this paper, we will not distinguish a skew partition $\lm$
with its diagram. Thus, the \emph{conjugate} of a skew partition
$\lm$, which we denote by $\lambda'/\mu'$, is defined to be the
transpose of the skew diagram $\lm$.

A \emph{(border) strip} or a \emph{ribbon} is a skew diagram with
an edgewise connected set of boxes that contains no $2\times 2$
block of boxes, where two boxes are said to be edgewise connected
if they share a common edge. One may impose a natural direction to
each strip. The starting box is always the one at the lower left
corner. Moreover, for each box that is not the last box in a
strip,  we  say that it goes right if the next box is to its
right, otherwise, we say that it goes up.

Recall that the \emph{content} of a box $(i,j)$ in a skew diagram
$\lm$ is given by $j-i$, and a \emph{diagonal} with content $c$ of
$\lm$ is the set of all the boxes in $\lm$ having content $c$.

Let $\lm$ be a skew diagram, a \emph{border strip decomposition}
of $\lm$ is a partitioning of the boxes of $\lm$ into pairwise
disjoint strips. Figure \ref{bd1} gives two examples of border
strip decompositions of the skew shape 6653/211. We now come to
the crucial definition of an outside decomposition of a skew shape
as introduced by Hamel and Goulden \cite{HG}. A border strip
decomposition of $\lm$ is said to be  an \emph{outside
decomposition} of $\lm$ if every strip in the decomposition has a
starting box on the left or bottom perimeter of the diagram and an
ending box on the right or top perimeter of the diagram, see
Figure \ref{bd1}.

\begin{figure}[h,t]
\begin{center}
\setlength{\unitlength}{20pt}
\begin{picture}(16,4)

\put(1,1){\line(1,0){2.1}} \put(1,1.7){\line(1,0){3.5}}
\put(1.7,2.4){\line(1,0){3.5}} \put(1.7,3.1){\line(1,0){3.5}}
\put(2.4,3.8){\line(1,0){2.8}}

\put(1,1){\line(0,1){0.7}}\put(1.7,1){\line(0,1){2.1}}
\put(2.4,1){\line(0,1){2.8}}\put(3.1,1){\line(0,1){2.8}}
\put(3.8,1.7){\line(0,1){2.1}}\put(4.5,1.7){\line(0,1){2.1}}
\put(5.2,2.4){\line(0,1){1.4}}

\multiput(1.35,1.35)(0.35,0){4}{\line(1,0){0.175}}
\multiput(2.75,3.45)(0.35,0){6}{\line(1,0){0.175}}
\multiput(4.15,2.75)(0.35,0){2}{\line(1,0){0.175}}

\multiput(2.75,1.35)(0,0.35){6}{\line(0,1){0.175}}
\multiput(2.05,1.875)(0,0.35){3}{\line(0,1){0.175}}
\multiput(3.45,1.875)(0,0.35){3}{\line(0,1){0.175}}
\multiput(4.15,1.875)(0,0.35){3}{\line(0,1){0.175}}

\put(10,1){\line(1,0){2.1}} \put(10,1.7){\line(1,0){3.5}}
\put(10.7,2.4){\line(1,0){3.5}} \put(10.7,3.1){\line(1,0){3.5}}
\put(11.4,3.8){\line(1,0){2.8}}

\put(10,1){\line(0,1){0.7}}\put(10.7,1){\line(0,1){2.1}}
\put(11.4,1){\line(0,1){2.8}}\put(12.1,1){\line(0,1){2.8}}
\put(12.8,1.7){\line(0,1){2.1}}\put(13.5,1.7){\line(0,1){2.1}}
\put(14.2,2.4){\line(0,1){1.4}}

\multiput(10.35,1.35)(0.35,0){2}{\line(1,0){0.175}}
\multiput(11.05,2.75)(0.35,0){2}{\line(1,0){0.175}}
\multiput(12.45,3.45)(0.35,0){4}{\line(1,0){0.175}}
\multiput(11.75,2.05)(0.35,0){2}{\line(1,0){0.175}}
\multiput(13.15,2.75)(0.35,0){2}{\line(1,0){0.175}}

\multiput(11.05,1.35)(0,0.35){4}{\line(0,1){0.175}}
\multiput(11.75,2.75)(0,0.35){3}{\line(0,1){0.175}}
\multiput(12.45,2.05)(0,0.35){4}{\line(0,1){0.175}}
\multiput(11.75,1.35)(0,0.35){2}{\line(0,1){0.175}}
\multiput(13.15,2.05)(0,0.35){2}{\line(0,1){0.175}}

\put(-0.5,0){a. A border strip decomposition} \put(9,0){b. An
outside decomposition}

\end{picture}
\end{center}
\caption{Border strip decompositions} \label{bd1}
\end{figure}

We note that there is slight difference between our description of
the outside decomposition and the original definition given by
Hamel and Goulden. Since the positions of the strips in an outside
decomposition are distinguishable, we may impose a canonical order
on the strips. However, in the original treatment of Hamel and
Goulden an order of the strips is taken into account and  they
have noted that any order plays the same role. As we will see
later in this paper, we may use the canonical order of the strips
by the contents of the ending boxes.

Given an outside decomposition $\Pi$ of a skew diagram $\lm$, and
a box $(i,j)$ in $\lm$, then there is a strip $\theta$ in $\Pi$
that contains $(i,j)$, and we say that the box $(i,j)$ \emph{goes
right} or \emph{goes up in $\lm$} if it goes right or goes up in
$\theta$. When $(i,j)$ is the ending box of $\theta$, then it is
on the top or the right perimeter of $\lm$, and we say that
$(i,j)$ \emph{goes up in $\lm$} if it is on the top perimeter of
$\lm$ (including the case when it is a corner of $\lm$), and
$(i,j)$ \emph{goes right in $\lm$} if it is on the right perimeter
of $\lm$. The directions of the boxes are important owing to the
following nested property.

\noindent {\bf Nested Property  (Hamel and Goulden \cite{HG}).}
The strips in any outside decomposition of $\lm$ are nested in the
sense that the boxes in the same diagonal of $\lm$ all go up or
all go right.

An immediate consequence of the above nested property is that the
contents of the ending(or starting) boxes of the strips in an
outside decomposition are different.  Hence we may use the
contents of the ending boxes to order the strips in an outside
decomposition. If there are $m$ strips in $\Pi$, we will denote
$\Pi$ by $\Pi=(\theta_1,\theta_2,\ldots,\theta_m)$, where
$\theta_1$ is the strip whose ending box has the largest content,
etc. In this notation, the outside decomposition of 6653/221 given
in Figure \ref{bd1}.b may be denoted by (4221/11,21,3322/211).

Let $\lm$ be a skew diagram. Recall that a \emph{block} of $\lm$
of a skew shape  is an edgewise connected component. If a skew
shape $\lm$ can be decomposed into $k$ blocks
$B_1,B_2,\ldots,B_k$, then an outside decomposition $\Pi$ of $\lm$
can be decomposed into $\Pi_1,\Pi_2,\ldots,\Pi_k$, where each
$\Pi_i$ is an outside decomposition of $B_i$, $i=1,2,\ldots,k$.
Therefore, without loss of generality one can only deal with
outside decompositions of edgewise connected skew shapes. In the
same vein, skew Schur function with disconnected shapes can be
expressed as a product of skew Schur functions with connected
shapes.

Based on the nested property, we are now ready to give
construction of the cutting strip of an outside decomposition of
an edgewise connected skew shape $\lm$.

\begin{deff} Let $\Pi$ be an outside decomposition of
an edgewise connected skew shape $\lm$. Suppoese that $\lm$ has
$d$ diagonals. The cutting strip $T$ of $\Pi$ is given by the
rule: for $i=1, 2, \ldots, d-1$,  the $i$-th box in $T$ goes right
or goes up according to whether the boxes in the $i$-th diagonal
of $\lm$ go right or go up with respect to $\Pi$.
\end{deff}

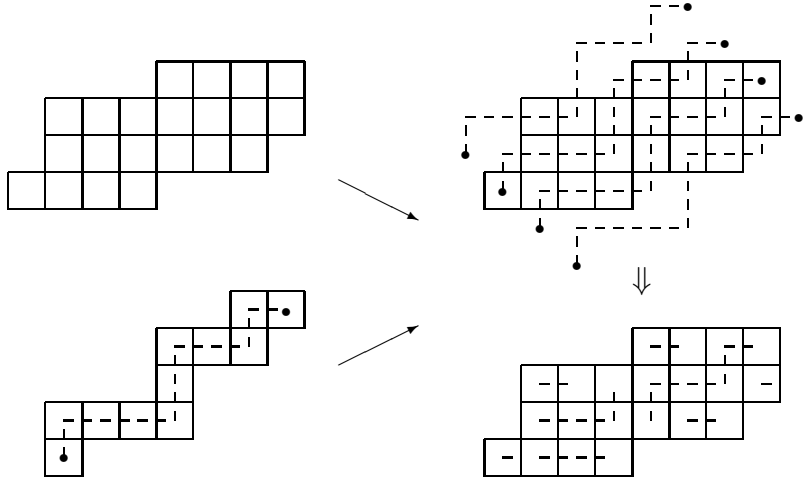
\begin{figure}[h,t]
\begin{center}
\setlength{\unitlength}{20pt}
\begin{picture}(16,9.6)

\put(0.75,5.95){\line(1,0){2.8}} \put(0.75,6.65){\line(1,0){4.9}}
\put(1.45,7.35){\line(1,0){4.9}} \put(1.45,8.05){\line(1,0){4.9}}
\put(3.55,8.75){\line(1,0){2.8}}

\put(0.75,5.95){\line(0,1){0.7}} \put(1.45,5.95){\line(0,1){2.1}}
\put(2.15,5.95){\line(0,1){2.1}} \put(2.85,5.95){\line(0,1){2.1}}
\put(3.55,5.95){\line(0,1){2.8}} \put(4.25,6.65){\line(0,1){2.1}}
\put(4.95,6.65){\line(0,1){2.1}} \put(5.65,6.65){\line(0,1){2.1}}
\put(6.35,7.35){\line(0,1){1.4}}

\put(1.45,0.9){\line(1,0){0.7}} \put(1.45,1.6){\line(1,0){2.8}}
\put(1.45,2.3){\line(1,0){2.8}} \put(3.55,3.0){\line(1,0){2.1}}
\put(3.55,3.7){\line(1,0){2.8}} \put(4.95,4.4){\line(1,0){1.4}}

\put(1.45,0.9){\line(0,1){1.4}} \put(2.15,0.9){\line(0,1){1.4}}
\put(2.85,1.6){\line(0,1){0.7}} \put(3.55,1.6){\line(0,1){2.1}}
\put(4.25,1.6){\line(0,1){2.1}} \put(4.95,3.0){\line(0,1){1.4}}
\put(5.65,3.0){\line(0,1){1.4}} \put(6.35,3.7){\line(0,1){0.7}}

\multiput(1.8,1.95)(0.35,0){6}{\line(1,0){0.175}}
\multiput(3.9,3.35)(0.35,0){4}{\line(1,0){0.175}}
\multiput(5.3,4.05)(0.35,0){2}{\line(1,0){0.175}}

\multiput(1.8,1.25)(0,0.35){2}{\line(0,1){0.175}}
\multiput(3.9,1.95)(0,0.35){4}{\line(0,1){0.175}}
\multiput(5.3,3.35)(0,0.35){2}{\line(0,1){0.175}}

\put(1.8,1.25){\circle*{0.125}} \put(6,4){\circle*{0.125}}

\put(7,6.5){\vector(2,-1){1.5}} \put(7,3){\vector(2,1){1.5}}
\put(12.55,4.45){$\Downarrow$}

\put(9.75,5.95){\line(1,0){2.8}} \put(9.75,6.65){\line(1,0){4.9}}
\put(10.45,7.35){\line(1,0){4.9}}
\put(10.45,8.05){\line(1,0){4.9}}
\put(12.55,8.75){\line(1,0){2.8}}

\put(9.75,5.95){\line(0,1){0.7}} \put(10.45,5.95){\line(0,1){2.1}}
\put(11.15,5.95){\line(0,1){2.1}}
\put(11.85,5.95){\line(0,1){2.1}}
\put(12.55,5.95){\line(0,1){2.8}}
\put(13.25,6.65){\line(0,1){2.1}}
\put(13.95,6.65){\line(0,1){2.1}}
\put(14.65,6.65){\line(0,1){2.1}}
\put(15.35,7.35){\line(0,1){1.4}}

\multiput(0.35,4.65)(0.7,-0.7){4}{\multiput(9.05,3.05)(0.35,0){6}{\line(1,0){0.175}}}
\multiput(0.35,4.65)(0.7,-0.7){4}{\multiput(11.15,4.45)(0.35,0){4}{\line(1,0){0.175}}}
\multiput(0.35,4.65)(0.7,-0.7){4}{\multiput(12.55,5.15)(0.35,0){2}{\line(1,0){0.175}}}

\multiput(0.35,4.65)(0.7,-0.7){4}{\multiput(9.05,2.35)(0,0.35){2}{\line(0,1){0.175}}}
\multiput(0.35,4.65)(0.7,-0.7){4}{\multiput(11.15,3.05)(0,0.35){4}{\line(0,1){0.175}}}
\multiput(0.35,4.65)(0.7,-0.7){4}{\multiput(12.55,4.45)(0,0.35){2}{\line(0,1){0.175}}}

\multiput(9.4,6.975)(0.7,-0.7){4}{\circle*{0.125}}
\multiput(13.6,9.775)(0.7,-0.7){4}{\circle*{0.125}}

\put(9.75,0.9){\line(1,0){2.8}} \put(9.75,1.6){\line(1,0){4.9}}
\put(10.45,2.3){\line(1,0){4.9}} \put(10.45,3){\line(1,0){4.9}}
\put(12.55,3.7){\line(1,0){2.8}}

\put(9.75,0.9){\line(0,1){0.7}} \put(10.45,0.9){\line(0,1){2.1}}
\put(11.15,0.9){\line(0,1){2.1}} \put(11.85,0.9){\line(0,1){2.1}}
\put(12.55,0.9){\line(0,1){2.8}} \put(13.25,1.6){\line(0,1){2.1}}
\put(13.95,1.6){\line(0,1){2.1}} \put(14.65,1.6){\line(0,1){2.1}}
\put(15.35,2.3){\line(0,1){1.4}}

\put(10.1,1.25){\line(1,0){0.175}}
\multiput(10.8,1.25)(0.35,0){4}{\line(1,0){0.175}}
\multiput(10.8,1.95)(0.35,0){4}{\line(1,0){0.175}}
\multiput(10.8,2.65)(0.35,0){2}{\line(1,0){0.175}}
\multiput(12.9,3.35)(0.35,0){2}{\line(1,0){0.175}}
\multiput(13.6,1.95)(0.35,0){2}{\line(1,0){0.175}}
\multiput(12.9,2.65)(0.35,0){4}{\line(1,0){0.175}}
\multiput(14.3,3.35)(0.35,0){2}{\line(1,0){0.175}}

\multiput(12.2,1.95)(0,0.35){2}{\line(0,1){0.175}}
\multiput(12.9,1.95)(0,0.35){2}{\line(0,1){0.175}}
\multiput(14.3,2.65)(0,0.35){2}{\line(0,1){0.175}}
\put(15,2.65){\line(1,0){0.175}}

\end{picture}
\end{center}
\caption{The cutting strip of an outside decomposition}
\label{od}
\end{figure}

Conversely, for any border strip with $d$ boxes, we may
reconstruct the outside decomposition $\Pi$ by peeling off $\lm$
while sliding the cutting strip along the diagonals. This
operation can be visualized by Figure \ref{od}.

Thus, we obtain the following correspondence.

\begin{theo}\label{bijection}
If $\lm$ is an edgewise connected skew shape with $d$ diagonals,
then there is a one-to-one correspondence between the outside
decompositions of $\lm$ and border strips with $d$ boxes.
\end{theo}

The above correspondence leads to the following formula for the
number of outside decompositions of an edgewise connected skew
shape.

\begin{coro}\label{enumod}
Let $\lm$ be an edgewise connected skew shape with $d$ diagonals.
Then $\lm$ has $2^{d-1}$ outside decompositions.
\end{coro}

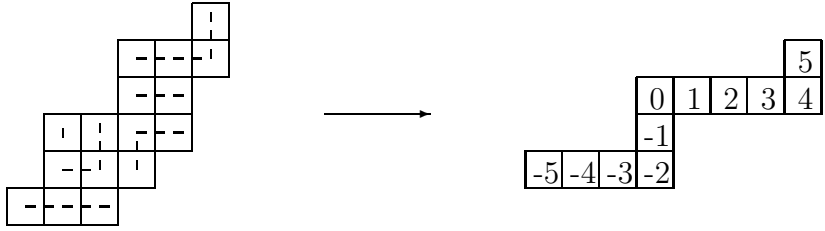
\begin{figure}[h,t]
\setlength{\unitlength}{20pt}
\begin{center}
\begin{picture}(15,5.5)(0,1)

\put(0,1){\line(1,0){2.1}} \put(0,1.7){\line(1,0){2.8}}
\put(0.7,2.4){\line(1,0){2.8}} \put(0.7,3.1){\line(1,0){2.8}}
\put(2.1,3.8){\line(1,0){2.1}} \put(2.1,4.5){\line(1,0){2.1}}
\put(3.5,5.2){\line(1,0){0.7}}

\put(0,1){\line(0,1){0.7}} \put(0.7,1){\line(0,1){2.1}}
\put(1.4,1){\line(0,1){2.1}} \put(2.1,1.0){\line(0,1){3.5}}
\put(2.8,1.7){\line(0,1){2.8}} \put(3.5,2.4){\line(0,1){2.8}}
\put(4.2,3.8){\line(0,1){1.4}}

\multiput(0.35,1.35)(0.35,0){5}{\line(1,0){0.175}}
\multiput(1.05,2.05)(0.35,0){2}{\line(1,0){0.175}}
\multiput(2.45,2.75)(0.35,0){3}{\line(1,0){0.175}}
\multiput(2.45,3.45)(0.35,0){3}{\line(1,0){0.175}}
\multiput(2.45,4.15)(0.35,0){4}{\line(1,0){0.175}}

\multiput(1.05,2.65)(0,0.35){1}{\line(0,1){0.175}}
\multiput(1.75,2.05)(0,0.35){3}{\line(0,1){0.175}}
\multiput(2.45,2.05)(0,0.35){2}{\line(0,1){0.175}}
\multiput(3.85,4.15)(0,0.35){3}{\line(0,1){0.175}}

\put(6.0,3.1){\vector(1,0){2.0}}

\put(9.8,1.7){\line(1,0){2.8}} \put(9.8,2.4){\line(1,0){2.8}}
\put(11.9,3.1){\line(1,0){3.5}} \put(11.9,3.8){\line(1,0){3.5}}
\put(14.7,4.5){\line(1,0){0.7}}

\put(9.8,1.7){\line(0,1){0.7}} \put(10.5,1.7){\line(0,1){0.7}}
\put(11.2,1.7){\line(0,1){0.7}} \put(11.9,1.7){\line(0,1){2.1}}
\put(12.6,1.7){\line(0,1){2.1}} \put(13.3,3.1){\line(0,1){0.7}}
\put(14.0,3.1){\line(0,1){0.7}} \put(14.7,3.1){\line(0,1){1.4}}
\put(15.4,3.1){\line(0,1){1.4}}

\put(9.95,1.8){-5} \put(10.65,1.8){-4} \put(11.35,1.8){-3}
\put(12.05,1.8){-2} \put(12.05,2.5){-1} \put(12.15,3.2){0}
\put(12.85,3.2){1} \put(13.55,3.2){2} \put(14.25,3.2){3}
\put(14.95,3.2){4} \put(14.95,3.9){5}

\end{picture}
\end{center}
\caption{Contents of a cutting strip} \label{d2t}
\end{figure}

The determinantal formula of Hamel and Goulden relies on the
operation $\#$ on border strips of an outside decomposition. We
will show that the $\#$ operation can be described in terms of the
cutting strip of the outside decomposition. To this end, we need
to record the contents of the diagonals in the construction of the
cutting strip. For example, the cutting strip in Figure \ref{d2t}
inherits the contents of the boxes in the original shape.

Given an outside decomposition $\Pi$ of a skew shape $\lm$ and a
strip $\theta$ in $\Pi$, let $\phi$ be the cutting strip of $\Pi$.
We denote the content of the starting box of $\theta$ and the
content of the ending box of $\theta$ respectively by $p(\theta)$
and $q(\theta)$, then $\theta$ forms a segment of the cutting
strip of $\Pi$ starting with the box with content $p(\theta)$ and
ending with the box with content $q(\theta)$, which is denoted by
$\phi[p(\theta), q(\theta)]$, or simply $[p(\theta), q(\theta)]$
if no confusion arises. Note that the contents of the cutting
strip are inherited from the contents in the original skew shape.
Furthermore, we may extend the notation $\phi[p,q]$ by the
following convention:
\begin{enumerate}
\item[(1)] if $p \leq q$, then $\phi[p,q]$ is a segment
of $\phi$ starting with the box having content $p$ and ending with
the box having content $q$;

\item[(2)] $\phi[q+1,q]=\varnothing$;

\item[(3)] if $p > q+1$, then $\phi[p,q]$ is
undefined.
\end{enumerate}

The difference between an empty strip and an undefined strip lies
in values of the corresponding Schur functions. For an empty
strip, the corresponding Schur function is defined to be 1, and
for the undefined strip, the Schur function takes the value 0.
Next we observe that the $\#$ operation on strips can be described
as extracting segments of the cutting strip. The proof of the
following fact is omitted since it is merely a technical
verification.

\begin{theo}\label{sharp}
Let $\lm$ be an edgewise connected skew shape, and $\Pi$ be an
outside decomposition of $\lm$. For any strips $\theta_i,\theta_j$
in $\Pi$, the strip $\theta_i\#\theta_j$, as defined in \cite{HG},
can be represented as $[p(\theta_j),q(\theta_i)]$.
\end{theo}

For example, the strips in Figure \ref{d2t} are $\theta_1=33/2$,
$\theta_2=2$, $\theta_3=21$, $\theta_4=22/1$, $\theta_5=1$,
$\theta_6=3$, and some of the strips obtained by the $\#$
operation are given below:
$$\begin{array}{lll}
\theta_1 \# \theta_5 = [-2,5] = 5511/4, & & \theta_5 \# \theta_1 = [2,-2] = {\rm undefined},\\
\theta_2 \# \theta_3 = [-1,2] = 31, & & \theta_3 \# \theta_2 = [1,1] = 1,\\
\theta_4 \# \theta_6 = [-5,-1] = 44/3, & & \theta_6 \# \theta_4 =
[-3,-3] = 1.
\end{array}$$

\begin{figure}[h,t]
\setlength{\unitlength}{20pt}
\begin{center}
\begin{picture}(12,2.5)(7,2)

\put(9.8,1.7){\line(1,0){2.8}} \put(9.8,2.4){\line(1,0){2.8}}
\put(11.9,3.1){\line(1,0){3.5}} \put(11.9,3.8){\line(1,0){3.5}}
\put(14.7,4.5){\line(1,0){0.7}}

\put(9.8,1.7){\line(0,1){0.7}} \put(10.5,1.7){\line(0,1){0.7}}
\put(11.2,1.7){\line(0,1){0.7}} \put(11.9,1.7){\line(0,1){2.1}}
\put(12.6,1.7){\line(0,1){2.1}} \put(13.3,3.1){\line(0,1){0.7}}
\put(14.0,3.1){\line(0,1){0.7}} \put(14.7,3.1){\line(0,1){1.4}}
\put(15.4,3.1){\line(0,1){1.4}}

\put(9.95,1.8){-5} \put(10.65,1.8){-4} \put(11.35,1.8){-3}
\put(12.05,1.8){\bf -2} \put(12.05,2.5){\bf -1}
\put(12.15,3.2){\bf 0} \put(12.85,3.2){\bf 1} \put(13.65,3.2){\bf
2} \put(14.25,3.2){\bf 3} \put(14.95,3.2){\bf
4}\put(14.95,3.9){\bf 5}

\end{picture}
\end{center}
\caption{$\theta_1 \# \theta_5$ as a segment $[-2,5]$ of the
cutting strip} \label{d2t2}
\end{figure}
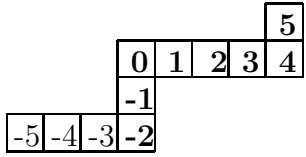

\section{Transformations of Determinants}

This section is concerned with a transformation theorem on the
determinantal formula of Hamel and Goulden \cite{HG}. In the proof
of the rim ribbon formula for Schur functions, Lascoux and Pragacz
give a transformation from the Giambelli determinant to the rim
ribbon determinant. Since the Hamel-Goulden determinant is a
unification of the four known determinantal identities, it is
natural to consider the transformation on the Hamel-Goulden
determinants, namely, from one outside decomposition to another.
From the cutting strip characterization of outside decompositions,
it follows that one can obtain any outside decomposition from
another by a sequence of basic transformations of changing the
directions of the boxes on a diagonal, which corresponds to the
operation of changing the direction of a box in the cutting strip.
This transformation is called {\em the twist transformation} on
border strips. Specifically, we use $\omega_i$ to denote the twist
transformation acting an outside decomposition  $\Pi$ by changing
the directions of the diagonal boxes with content $i$.

We may restrict our attention to the determinantal transformations
that correspond to the twist transformations on the cutting
strips, or equivalently on the outside decompositions. Let us
review the Hamel-Goulden determinantal formula. Without loss of
generality, we may assume that $\lm$ is an edgewise connected skew
shape with $d$ diagonals. Let
$\Pi=(\theta_1,\theta_2,\ldots,\theta_m)$ be an outside
decomposition of $\lm$. Then we have

\begin{theo}[Hamel and Goulden \cite{HG}]
The skew Schur function $s_\lm(X)$ can be evaluated by the
following determinant:
\begin{equation}\label{eq1}
D(\Pi)=
\det\left[\begin{matrix}s_{\theta_1\#\theta_1}&s_{\theta_1\#\theta_2}&\ldots&
s_{\theta_1\#\theta_m}\\
s_{\theta_2\#\theta_1}&s_{\theta_2\#\theta_2}&\ldots&s_{\theta_2\#\theta_m}\\
\vdots&\vdots&\ldots&\vdots\\
s_{\theta_m\#\theta_1}&s_{\theta_m\#\theta_2}&\ldots&s_{\theta_m\#\theta_m}
\end{matrix}\right],
\end{equation}
where $s_\varnothing=1$ and $s_{\theta_i \# \theta_j}=0$ if
$\theta_i \# \theta_j$ is undefined.
\end{theo}

We recall that the above determinantal formula includes the
Jacobi-Trudi identity, the dual Jacobli-Trudi identity (the
N$\rm{\ddot{a}}$gelsbach-Kostka identity), the Giambelli identity
and the rim ribbon identity as special cases. Note that for any
order of the strips in the outside decomposition, the above
determinant remains invariant. The corresponding cutting strips
for the above special cases are the horizontal strip, the vertical
strip, the maximal hook, and the maximal outside rim ribbon, see
\cite{HG}. As noticed by Lascoux and Pragacz \cite{LP2}, one may
obtain the rim ribbon identity from the Giambelli identity, and
vice versa. We may expect that the above mentioned four identities
can be transformed from one to another. In fact, this goal can be
achieved in a more general setting in terms of the transformation
theorem on the Hamel-Goulden identity.

For an outside decomposition $\Pi=(\theta_1, \theta_2, \ldots,
\theta_m)$ of $\lm$, we define its inversion number by
\[\inv(\Pi)=|\{(i,j); p(\theta_i) > p(\theta_j),\ q(\theta_i) < q(\theta_j)\}|.\]
By the cutting strip characterization of the $\#$ operation, we
may rewrite (\ref{eq1}) as
\begin{equation}\label{eq2}
(-1)^{\inv(\Pi)} \det\left[
\begin{matrix}s_{[p_1,\ q_1]}&s_{[p_2,\ q_1]}&\ldots&
s_{[p_m,\ q_1]}\\
s_{[p_1,\ q_2]}&s_{[p_2,\ q_2]}&\ldots&s_{[p_m,\ q_2]}\\
\vdots&\vdots&\ldots&\vdots\\
s_{[p_1,\ q_m]}&s_{[p_2,\ q_m]}&\ldots&s_{[p_m,\ q_m]}
\end{matrix}\right].
\end{equation}
 Where $p_1>p_2>\ldots>p_m$ are the rearrangement
of $p(\theta_1), \ldots, p(\theta_m)$ in decreasing order and
$q_i=q(\theta_i)$ for $1 \leq i \leq m$. Since
$p_1>p_2>\ldots>p_m$ and $q_1>q_2>\ldots>q_m$, it is clear that if
$s_{[p_i,q_j]}=0$ then $s_{[p_i,\ q_{j\,'}]}=0$ and
$s_{[p_{i\,'},\ q_j]}=0$ for $j\leq j\,'\leq m$ and $1\leq
i\,'\leq i$. The above canonical form is more convenient for the
construction of  the determinantal transformations corresponding
to the twist transformations on outside decompositions.

Now the question becomes how to transform the determinant
(\ref{eq2}) for an outside decomposition $\Pi$ to the determinant
corresponding to the outside decomposition $\omega_i(\Pi)$. To
accomplish this goal, we need an important product rule for the
skew Schur functions due to Zelevinsky \cite{Z}, which follows
from the \emph{Jeu de Taquin} of Sch\"utzenberger \cite{Sch}, and
can be verified by definition. This product rule is also used by
Lascoux and Pragacz \cite{LP2} in their transformation from the
Giambelli identity to the rim ribbon identity.

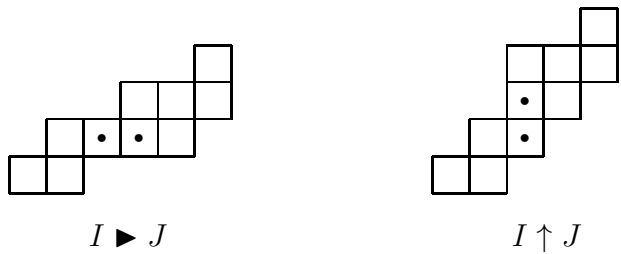
\begin{figure}[h,t]
\begin{center}
\setlength{\unitlength}{20pt}
\begin{picture}(14,4.5)

\put(1,1){\line(1,0){1.4}} \put(1,1.7){\line(1,0){3.5}}
\put(1.7,2.4){\line(1,0){3.5}} \put(3.1,3.1){\line(1,0){2.1}}
\put(4.5,3.8){\line(1,0){0.7}}

\put(1,1){\line(0,1){0.7}}\put(1.7,1){\line(0,1){1.4}}
\put(2.4,1){\line(0,1){1.4}}\put(3.1,1.7){\line(0,1){1.4}}
\put(3.8,1.7){\line(0,1){1.4}}
\put(4.5,1.7){\line(0,1){2.1}}\put(5.2,2.4){\line(0,1){1.4}}

\put(9,1){\line(1,0){1.4}} \put(9,1.7){\line(1,0){2.1}}
\put(9.7,2.4){\line(1,0){2.1}} \put(10.4,3.1){\line(1,0){2.1}}
\put(10.4,3.8){\line(1,0){2.1}} \put(11.8,4.5){\line(1,0){0.7}}

\put(9,1){\line(0,1){0.7}}\put(9.7,1){\line(0,1){1.4}}
\put(10.4,1){\line(0,1){2.8}}\put(11.1,1.7){\line(0,1){2.1}}
\put(11.8,2.4){\line(0,1){2.1}}\put(12.5,3.1){\line(0,1){1.4}}

\put(2.75,2.05){\circle*{0.125}} \put(3.45,2.05){\circle*{0.125}}
\put(10.75,2.05){\circle*{0.125}}
\put(10.75,2.75){\circle*{0.125}}

\put(2.5,0){$I\blacktriangleright J$ } \put(10.5,0){$I\uparrow J$}
\end{picture}
\end{center}
\caption{The $\blacktriangleright$ and $\uparrow$ operator for
skew diagrams}\label{bu}
\end{figure}

Given two skew diagrams $I$ and $J$, $I\blacktriangleright J$ is
the diagram obtained by gluing the lower left-hand corner box of
diagram $J$ to the right of the upper right-hand corner box of
diagram $I$, $I\uparrow J$ is the diagram obtained by gluing the
lower left-hand corner box of diagram $J$ to the top of the upper
right-hand corner box of diagram $I$. For example, if $I=32/1$ and
$J=332/2$, then $I\blacktriangleright J$ and $I\uparrow J$ are
given by Figure \ref{bu}

\begin{lemm}\label{llas}
Let $I$ and $J$ be two skew diagrams. Then
$$s_{I}s_{J}=s_{I \blacktriangleright J}+s_{I \uparrow J}.$$
\end{lemm}

As far as this paper is concerned, we will only need the above
property for ribbon Schur functions, or the strip Schur functions.
More specifically, we need to multiply a row of the Hamel-Goulden
determinant by a ribbon Schur function and subtract the product
from another row to obtain a new ribbon Schur function.

For two strips $\alpha$ and $\beta$, one has two ways to combine
them into another ribbon, one is $\alpha \blacktriangleright
\beta$ and the other is $\alpha \uparrow \beta$. It turns out that
the basic transformation on the cutting strip is an interchange of
the operations $\alpha \blacktriangleright \beta$ and $\alpha
\uparrow \beta$. This is the reason why we call $\omega_i$ the
twist transformation.

\begin{theo}\label{thmmain}
If $\lm$ is an edgewise connected skew shape with $d$ diagonals,
let $\Pi$ be an outside decomposition of $\lm$, then for each
diagonal of $\lm$ with content $i$, we can obtain
$D(\omega_i(\Pi))$ from $D(\Pi)$ by elementary determinantal
operations.
\end{theo}

We remark that the elementary determinantal operations mentioned
in the above theorem actually involve three types, i.e.,
interchanging two rows or columns; multiplying a row or a column
by a non-zero scalar, or adding a multiple of one row or one
column to another row or column; and from two determinants, one
can construct a bigger determinant by forming a bigger diagonal
block matrix, or one can reduce a diagonal determinant into two
smaller determinants, one of which may be the identity matrix.

\begin{figure}[h,t]
\setlength{\unitlength}{20pt}
\begin{center}
\begin{picture}(15,12)

\put(3.8,7){\line(1,0){0.7}} \put(3.1,7.7){\line(1,0){1.4}}
\put(3.1,8.4){\line(1,0){1.4}} \put(3.1,9.1){\line(1,0){0.7}}
\put(1.7,9.1){\line(1,0){0.7}} \put(1,9.8){\line(1,0){1.4}}
\put(1,10.5){\line(1,0){1.4}} \put(1,11.2){\line(1,0){0.7}}

\put(0.4,9.9){L} \put(8.4,9.9){L}

\put(1,9.8){\line(0,1){1.4}} \put(1.7,9.1){\line(0,1){2.1}}
\put(2.4,9.1){\line(0,1){1.4}} \put(3.1,7.7){\line(0,1){1.4}}
\put(3.8,7){\line(0,1){2.1}} \put(4.5,7){\line(0,1){1.4}}

\put(2.5,9.2){$\ddots$}\put(2.5,8.5){$\ddots$}

\put(1.35,10.15){\circle*{0.125}} \put(4.15,7.35){\circle*{0.125}}
\put(3.45,8.05){\circle*{0.125}} \put(2.05,9.45){\circle*{0.125}}

\put(11.8,7){\line(1,0){1.4}} \put(11.1,7.7){\line(1,0){2.1}}
\put(11.1,8.4){\line(1,0){1.4}} \put(9.7,9.1){\line(1,0){1.4}}
\put(9,9.8){\line(1,0){2.1}} \put(9,10.5){\line(1,0){1.4}}

\put(9,9.8){\line(0,1){0.7}} \put(9.7,9.1){\line(0,1){1.4}}
\put(10.4,9.1){\line(0,1){1.4}} \put(11.1,7.7){\line(0,1){0.7}}
\put(11.1,9.1){\line(0,1){0.7}} \put(11.8,7){\line(0,1){1.4}}
\put(12.5,7){\line(0,1){1.4}} \put(13.2,7){\line(0,1){0.7}}

\put(11.2,8.5){$\ddots$} \put(10.5,8.5){$\ddots$}

\put(9.35,10.15){\circle*{0.125}}
\put(12.15,7.35){\circle*{0.125}}
\put(11.45,8.05){\circle*{0.125}}
\put(10.05,9.45){\circle*{0.125}}

\put(2.5,6){Type 1 } \put(10.5,6){Type 2}

\put(0.4,3.9){L} \put(9.1,3.2){L}

\put(3.8,1){\line(1,0){0.7}} \put(3.1,1.7){\line(1,0){1.4}}
\put(3.1,2.4){\line(1,0){1.4}} \put(3.1,3.1){\line(1,0){0.7}}
\put(1.7,3.1){\line(1,0){0.7}} \put(1,3.8){\line(1,0){1.4}}
\put(1,4.5){\line(1,0){1.4}}

\put(1,3.8){\line(0,1){0.7}} \put(1.7,3.1){\line(0,1){1.4}}
\put(2.4,3.1){\line(0,1){1.4}} \put(3.1,1.7){\line(0,1){1.4}}
\put(3.8,1){\line(0,1){2.1}} \put(4.5,1){\line(0,1){1.4}}

\put(2.5,3.2){$\ddots$}\put(2.5,2.5){$\ddots$}

\put(1.35,4.15){\circle*{0.125}} \put(4.15,1.35){\circle*{0.125}}
\put(3.45,2.05){\circle*{0.125}} \put(2.05,3.45){\circle*{0.125}}

\put(11.8,1){\line(1,0){1.4}} \put(11.1,1.7){\line(1,0){2.1}}
\put(11.1,2.4){\line(1,0){1.4}} \put(9.7,3.1){\line(1,0){1.4}}
\put(9.7,3.8){\line(1,0){1.4}} \put(9.7,4.5){\line(1,0){0.7}}

\put(9.7,3.1){\line(0,1){1.4}} \put(10.4,3.1){\line(0,1){1.4}}
\put(11.1,1.7){\line(0,1){0.7}} \put(11.1,3.1){\line(0,1){0.7}}
\put(11.8,1){\line(0,1){1.4}} \put(12.5,1){\line(0,1){1.4}}
\put(13.2,1){\line(0,1){0.7}}

\put(11.2,2.5){$\ddots$} \put(10.5,2.5){$\ddots$}

\put(12.15,1.35){\circle*{0.125}}
\put(11.45,2.05){\circle*{0.125}}
\put(10.05,3.45){\circle*{0.125}}

\put(2.5,0){Type 3 } \put(10.5,0){Type 4}
\end{picture}
\end{center}
\caption{Four possible types of diagonals of $\lm$} \label{4types}
\end{figure}
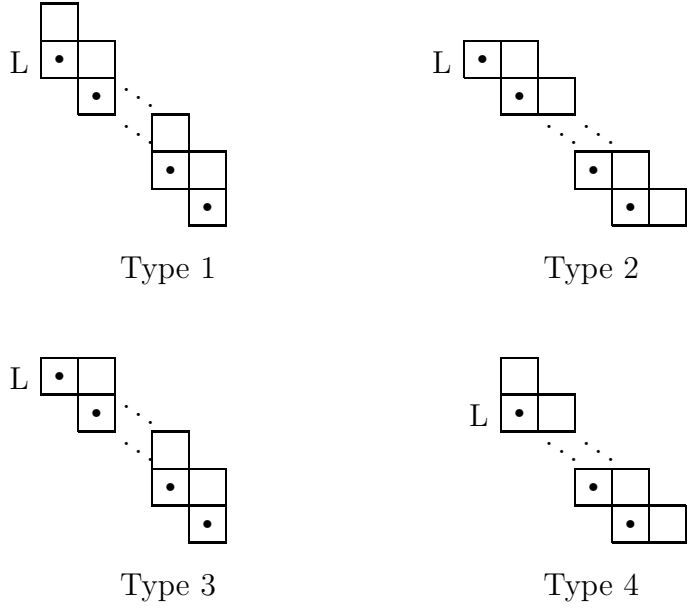

Note that if $L$ is a diagonal of an edgewise connected skew
diagram $\lm$, then there are four possible diagonal types
classified by whether the first diagonal box has a box at the top,
and whether the last diagonal box has a box on the right. These
four types are depicted  by Figure \ref{4types}.

We now proceed to give the proof of Theorem \ref{thmmain}. There
are essentially two different types (Type 1 and Type 3), because
Type 2 is similar to Type 1, and Type 4 is similar to Type 3. For
Type 1  and Type 2 diagonals, the determinantal transformations
involve changing a determinant to another determinant of different
order. We need the following lemma for counting the number of
inversions of a sequence of pairs of numbers, which will be used
to deal with the sign change of determinants after permuting rows
and columns. Given a sequence of $m$ pairs
$P=\{(a_1,b_1),(a_2,b_2),\ldots,(a_m,b_m)\}$ with
$a_1,a_2,\ldots,a_m$ being distinct and $b_1,b_2,\ldots,b_m$ being
distinct. We define inversion number of the sequence by
\[\inv(P)=|\{(i,j); a_i>a_j,\ b_i<b_j\}|.\]
It can be easily seen that the inversion depends only on the set
of pairs in a sequence. We choose the sequence notation for the
sake of presentation.

\begin{lemm}\label{inv}
Let $P=\{(a_1,b_1),(a_2,b_2),\ldots,(a_m,b_m)\}$ be a sequence of
$m$ pairs such that $a_1, a_2, \ldots, a_m$ are distinct and $b_1,
b_2, \ldots, b_m$ are distinct. Let $\sigma$ be a permutation of
$1,2,\ldots,m$ acting on the indices of $a_1, a_2, \ldots, a_m$,
$\sigma(a_i)=a_{\sigma(i)}$, then we have
\[\inv(\{(\sigma(a_1),b_1),(\sigma(a_2),b_2),\ldots,(\sigma(a_m),b_m)\})\equiv \inv(\sigma)+\inv(P)
\ \textnormal{(mod 2)}.\]
\end{lemm}

\noindent \proof Suppose $b_{i_1}<b_{i_2}<\ldots<b_{i_m}$ is the
reordering of $b_1,b_2,\ldots,b_m$, then we have
\[ \inv(P)=\inv((a_{i_1},a_{i_2},\ldots,a_{i_m}))=|\{(k,l);a_{i_k}>a_{i_l},
k<l\}|.\]
and
\begin{eqnarray*}
\lefteqn{\inv(\{(\sigma(a_1),b_1),(\sigma(a_2),b_2),\ldots,(\sigma(a_m),b_m)\})}\\
&&=\inv((a_{\sigma(i_1)},a_{\sigma(i_2)},\ldots,a_{\sigma(i_m)}))\\
&&\equiv \inv(\sigma)+\inv((a_{i_1},a_{i_2},\ldots,a_{i_n}))\
\textnormal{(mod 2)}\\
&&=\inv(\sigma)+\inv(P).
\end{eqnarray*}
This completes the proof. \qed

We are now ready to present the proof of Theorem \ref{thmmain} for
Type 1 diagonals.  Let $\Pi=(\theta_1,\theta_2,\ldots,\theta_m)$
be an outside decomposition of an edgewise connected skew diagram
$\lm$, and let $\phi$ be the cutting strip of $\Pi$. Recall that
the sequence of the contents of the starting boxes of the strips
in $\Pi$ is $p_1>p_2>\ldots>p_m$ and the sequence of the contents
of the ending boxes of the strips in $\Pi$ is
$q_1>q_2>\ldots>q_m$.

Suppose that $\lm$ has a Type 1 diagonal $L$ with content $i$, and
suppose that there are $r$ boxes in $L$. We may assume without
loss of generality that the diagonal boxes all go up, since we may
reverse the transformation process for the case that the diagonal
boxes all go right. Thus the cutting strip of $\Pi$ may be written
as $[p_m,i]\uparrow[i+1,q_1]$ and the cutting strip of
$\omega_i(\Pi)$ may be written as
$[p_m,i]\blacktriangleright[i+1,q_1]$. Notice that the twist
transformation $\omega_i$ does not change  the contents of the
cutting strip of $\Pi$.

\begin{figure}[h,t]
\begin{center}
\setlength{\unitlength}{20pt}
\begin{picture}(14,6.5)

\put(3.8,1){\line(1,0){0.7}}\put(3.1,1.7){\line(1,0){1.4}}
\put(3.1,2.4){\line(1,0){1.4}}\put(3.1,3.1){\line(1,0){0.7}}
\put(1.7,3.1){\line(1,0){0.7}}\put(1,3.8){\line(1,0){1.4}}
\put(1,4.5){\line(1,0){1.4}}\put(1,5.2){\line(1,0){0.7}}

\put(1,3.8){\line(0,1){1.4}}\put(1.7,3.1){\line(0,1){2.1}}
\put(2.4,3.1){\line(0,1){1.4}}\put(3.1,1.7){\line(0,1){1.4}}
\put(3.8,1){\line(0,1){2.1}}\put(4.5,1){\line(0,1){1.4}}
\thicklines
\put(4.15,1.35){\vector(0,1){0.7}}\put(3.45,2.05){\vector(0,1){0.7}}
\put(2.05,3.45){\vector(0,1){0.7}}\put(1.35,4.15){\vector(0,1){0.7}}
\thinlines
\put(4.15,2.05){\line(0,1){0.7}}\put(3.45,1.35){\line(1,0){0.7}}
\put(3.45,2.75){\line(0,1){0.7}}\put(2.75,2.05){\line(1,0){0.7}}
\put(2.05,4.15){\line(0,1){0.7}}\put(1.35,3.45){\line(1,0){0.7}}
\put(1.35,4.85){\line(0,1){0.7}}\put(0.65,4.15){\line(1,0){0.7}}

\put(4.15,2.85){\small$\theta_{i_r}$}\put(3.45,3.55){\small$\theta_{i_{r-1}}$}
\put(2.05,4.95){\small$\theta_{i_2}$}\put(1.35,5.65){\small$\theta_{i_1}$}
\put(2.75,1.15){\small$\theta_{i_r}$}\put(1.65,1.85){\small$\theta_{i_{r-1}}$}
\put(0.65,3.25){\small$\theta_{i_2}$}\put(-0.05,3.95){\small$\theta_{i_1}$}

\put(4.85,1.0){\small $i$} \put(-0.25,4.75){\small $i+1$}

\put(2.5,3.2){$\ddots$}\put(2.5,2.5){$\ddots$}

\put(11.8,1){\line(1,0){0.7}}\put(11.1,1.7){\line(1,0){1.4}}
\put(11.1,2.4){\line(1,0){1.4}}
\put(9.7,3.1){\line(1,0){1.4}}\put(9,3.8){\line(1,0){2.1}}
\put(9,4.5){\line(1,0){1.4}}\put(9,5.2){\line(1,0){0.7}}

\put(9,3.8){\line(0,1){1.4}}\put(9.7,3.1){\line(0,1){2.1}}
\put(10.4,3.1){\line(0,1){1.4}}\put(11.1,1.7){\line(0,1){0.7}}
\put(11.1,3.1){\line(0,1){0.7}}
\put(11.8,1){\line(0,1){1.4}}\put(12.5,1){\line(0,1){1.4}}

\thicklines \put(11.45,2.05){\vector(1,0){0.7}}
\put(10.05,3.45){\vector(1,0){0.7}}\put(9.35,4.15){\vector(1,0){0.7}}
\thinlines
\put(12.15,2.05){\line(0,1){0.7}}\put(10.75,2.05){\line(1,0){0.7}}
\put(10.75,3.45){\line(0,1){0.7}}\put(9.35,3.45){\line(1,0){0.7}}
\put(10.05,4.15){\line(0,1){0.7}}\put(8.65,4.15){\line(1,0){0.7}}

\put(12.15,2.85){\small$\theta_{i_r}$}\put(10.75,4.25){\small$\theta_{i_3}$}
\put(10.05,4.95){\small$\theta_{i_2}$}\put(9.65,1.85){\small$\theta_{i_{r-1}}$}
\put(8.65,3.25){\small$\theta_{i_2}$}\put(7.95,3.95){\small$\theta_{i_1}$}
\put(10.75,1.15){\small$\theta_{i_r}$}

\put(11.2,2.5){$\ddots$}\put(10.5,2.5){$\ddots$}

\put(9.35,4.85){\circle*{0.125}} \put(12.15,1.35){\circle*{0.125}}
\put(9.35,4.85){\line(0,1){0.7}} \put(11.45,1.35){\line(1,0){0.7}}
\put(9.35,5.65){\small$\theta_{i_1}$}

\put(12.85,1.0){\small $i$} \put(7.75,4.75){\small $i+1$}

\put(2.5,0){$L$ } \put(10.5,0){$\omega_i(L)$}
\end{picture}
\end{center}
\caption{$\omega_i$ acts on a Type 1 diagonal $L$ }\label{D1}
\end{figure}

Since $\omega_i$ only changes the strips which contain a box in
$L$, we may suppose that $\theta_{i_t}$, $1 \leq t \leq r$, is the
strip in $\Pi$ that contains the $t$-th diagonal box in $L$. Then
we have $p(\theta_{i_t})\leq i < q(\theta_{i_t})$, for $1\leq t
\leq r$. As illustrated in Figure \ref{D1}, under the operation of
$\omega_i$, the strip
\[ \theta_{i_1} =
[p(\theta_{i_1}),q(\theta_{i_1})]=[p(\theta_{i_1}),i]\uparrow[i+1,q(\theta_{i_1})]\]
 breaks into two strips
 \[
[p(\theta_{i_1}),q(\theta_{i_2})]=[p(\theta_{i_1}),i]\blacktriangleright[i+1,q(\theta_{i_2})]
\quad \mbox{ and} \quad [i+1,q(\theta_{i_1})].\]
 If $r > 1$ then the last
strip
\[ \theta_{i_r} =
[p(\theta_{i_r}),q(\theta_{i_r})]=[p(\theta_{i_r}),i]\uparrow[i+1,q(\theta_{i_r})]\]
will be cut off into $[p(\theta_{i_r}),i]$, and the other strips
\[ \theta_{i_t} =
[p(\theta_{i_t}),q(\theta_{i_t})]=[p(\theta_{i_t}),i]\uparrow[i+1,q(\theta_{i_t})],\
2 \leq t \leq r-1, \] will be twisted into
\[[p(\theta_{i_t}),q(\theta_{i_{t+1}})]=[p(\theta_{i_t}),i]\blacktriangleright[i+1,q(\theta_{i_{t+1}})].\]
By Lemma \ref{inv} one can verify that $\inv(\omega_i(\Pi)) \equiv
\inv(\Pi)\ (\textnormal{mod 2})$.

Since $p(\theta_{i_t})< i+1$ for $1 \leq t \leq r$, there is a
unique integer $k,\ 1 \leq k \leq m-r+1$, such that  $p_k < i+1$
and $p_{k-1}
> i+1$ where $p_0$ is defined to be $\infty$, if necessary.
Therefore, there are exactly $k-1$ strips in $\Pi$ for which the
contents of starting boxes are bigger than $i+1$, and there are
exactly $k-1+r$ strips in $\Pi$ for which the contents of ending
boxes are bigger than $i$. It follows that $q_{k+r-1}>i$ and
$q_{k+r}<i$ if $k+r<m+1$. Moreover, the sequence of the contents
of the starting boxes of the strips in $\omega_i(\Pi)$ becomes
$p_1>\ldots>p_{k-1}>i+1>p_k>\ldots>p_m$, and the sequence of the
contents of the ending boxes of the strips in $\omega_i(\Pi)$
becomes $q_1>\ldots>q_{k+r-1}>i>q_{k+r}>\ldots>q_m$. As noted
before, the determinant $D(\Pi)$ has the following canonical form:
\begin{equation}\label{H1}
(-1)^{\inv(\Pi)} \det\left[
\begin{matrix}
A&B\\
0&C
\end{matrix}\right],
\end{equation}
where
\begin{equation*}
A=\left[\begin{matrix} s_{[p_1,\ q_1]}&\ldots&
s_{[p_{k-1},\ q_1]}\\
\vdots&\ldots&\vdots\\
s_{[p_1,\ q_{k+r-1}]}&\ldots& s_{[p_{k-1},\ q_{k+r-1}]}
\end{matrix}\right]_{(k+r-1)\times (k-1)},
\end{equation*}

\begin{equation*}
B=\left[\begin{matrix}s_{[p_k,\ i]\uparrow[i+1,\ q_1]}&\ldots&
s_{[p_m,\ i]\uparrow[i+1,\ q_1]}\\
\vdots&\ldots&\vdots\\
s_{[p_k,\ i]\uparrow[i+1,\ q_{k+r-1}]}&\ldots& s_{[p_m,\
i]\uparrow[i+1,\ q_{k+r-1}]}
\end{matrix}\right]_{(k+r-1)\times (m-k+1)},
\end{equation*}

\begin{equation*}
C=\left[\begin{matrix}s_{[p_k,\ q_{k+r}]}&\ldots&
s_{[p_m,\ q_{k+r}]}\\
\vdots&\ldots&\vdots\\
s_{[p_k,\ q_m]}&\ldots&s_{[p_m,\ q_m]}
\end{matrix}\right]_{(m-k-r+1)\times (m-k+1)} .
\end{equation*}

Using the sequences of the contents of starting boxes and ending
boxes of the strips of $\omega_i(\Pi)$, we may obtain the
following canonical form of  $D(\omega_i(\Pi))$:

\begin{equation}\label{H2}
(-1)^{\inv(\omega_i(\Pi))}\det\left[
\begin{matrix}
\text{\huge \emph{A}}&\begin{array}{c}s_{[i+1,\ q_1]}\\ \vdots\\
s_{[i+1,\ q_{k+r-1}]}\end{array}&
\text{\huge \emph{B}}'\\
\begin{array}{ccc}0&\ldots&0\end{array}&
1&\begin{array}{ccc}s_{[p_k,\ i]}&\ldots&s_{[p_m,\ i]}\end{array}\\
\text{\huge \emph{0}}&\begin{array}{c}0\\ \vdots\\
0\end{array}& \text{\huge \emph{C}}
\end{matrix}\right],
\end{equation}
where
\begin{equation*}
B'=\left[\begin{matrix}s_{[p_k,\ i]\blacktriangleright[i+1,\
q_1]}&\ldots&
s_{[p_m,\ i]\blacktriangleright[i+1,\ q_1]}\\
\vdots&\ldots&\vdots\\
s_{[p_k,\ i]\blacktriangleright[i+1,\ q_{k+r-1}]}&\ldots&
s_{[p_m,\ i]\blacktriangleright[i+1,\ q_{k+r-1}]}
\end{matrix}\right]_{(k+r-1)\times (m-k+1)}.
\end{equation*}

Now we need to transform \eqref{H1} into the following determinant
by adding a row and a column as shown below. This operation
involves a sign $(-1)^{2k+r}$:
\begin{equation}\label{H21}
(-1)^{\inv(\Pi)}\det\left[
\begin{matrix}
\text{\huge \emph{A}}&\begin{array}{c}0\\ \vdots\\
0\end{array}&
\text{\huge \emph{B}}\\
\begin{array}{ccc}0&\ldots&0\end{array}&
1&\begin{array}{ccc}s_{[p_k,\ i]}&\ldots&s_{[p_m,\ i]}\end{array}\\
\text{\huge \emph{0}}&\begin{array}{c}0\\ \vdots\\
0\end{array}& \text{\huge \emph{C}}
\end{matrix}\right]
\end{equation}
Applying Lemma \ref{llas}, we have
\[s_{[p_t,\ i]\blacktriangleright[i+1,\ q_l]}=s_{[p_t,\ i]}s_{[i+1,\ q_l]}-s_{[p_t,\ i]\uparrow[i+1,\ q_l]}\]
for $k \leq t \leq m$ and $1 \leq l \leq k+r-1$. The above
identity enables us to use the $(k+r)$-th row of \eqref{H21} to
twist the entries in $B$ by elementary determinantal operations.
Such operations will change $B$ into $-B'$. Multiplying the last
$m-k+2$ columns and the last $m-k-r+2$ rows of \eqref{H21} by
$-1$, we may derive  \eqref{H2} from \eqref{H21}. It follows that
we can transform $D(\Pi)$ into $D(\omega_i(\Pi))$ by elementary
determinantal operations. This completes the proof of Theorem
\ref{thmmain} for Type 1 diagonals.

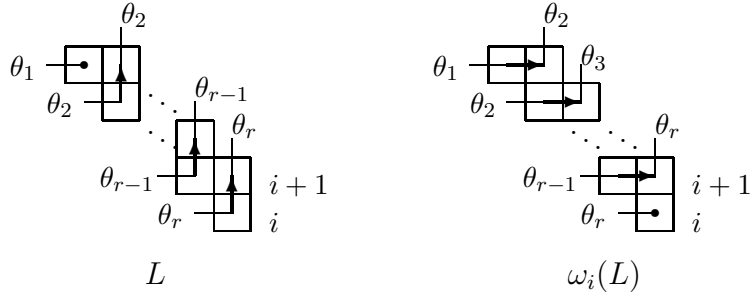
\begin{figure}[h,t]
\begin{center}
\setlength{\unitlength}{20pt}
\begin{picture}(14,6)
\put(3.8,1){\line(1,0){0.7}}\put(3.1,1.7){\line(1,0){1.4}}
\put(3.1,2.4){\line(1,0){1.4}}\put(3.1,3.1){\line(1,0){0.7}}
\put(1.7,3.1){\line(1,0){0.7}}\put(1,3.8){\line(1,0){1.4}}
\put(1,4.5){\line(1,0){1.4}}

\put(1,3.8){\line(0,1){0.7}}\put(1.7,3.1){\line(0,1){1.4}}
\put(2.4,3.1){\line(0,1){1.4}}\put(3.1,1.7){\line(0,1){1.4}}
\put(3.8,1){\line(0,1){2.1}}\put(4.5,1){\line(0,1){1.4}}

\thicklines
\put(4.15,1.35){\vector(0,1){0.7}}\put(3.45,2.05){\vector(0,1){0.7}}
\put(2.05,3.45){\vector(0,1){0.7}} \thinlines

\put(4.15,2.05){\line(0,1){0.7}}\put(3.45,1.35){\line(1,0){0.7}}
\put(3.45,2.75){\line(0,1){0.7}}\put(2.75,2.05){\line(1,0){0.7}}
\put(2.05,4.15){\line(0,1){0.7}}\put(1.35,3.45){\line(1,0){0.7}}

\put(4.15,2.85){\small$\theta_r$}\put(3.45,3.55){\small$\theta_{r-1}$}
\put(2.05,4.95){\small$\theta_2$}\put(2.75,1.15){\small$\theta_r$}
\put(1.65,1.85){\small$\theta_{r-1}$}\put(0.65,3.25){\small$\theta_2$}
\put(-0.05,3.95){\small$\theta_1$}

\put(4.85,1.0){\small $i$} \put(4.85,1.7){\small $i+1$}

\put(2.5,3.2){$\ddots$}\put(2.5,2.5){$\ddots$}

\put(11.8,1){\line(1,0){0.7}}\put(11.1,1.7){\line(1,0){1.4}}
\put(11.1,2.4){\line(1,0){1.4}}
\put(9.7,3.1){\line(1,0){1.4}}\put(9,3.8){\line(1,0){2.1}}
\put(9,4.5){\line(1,0){1.4}}

\put(9,3.8){\line(0,1){0.7}}\put(9.7,3.1){\line(0,1){1.4}}
\put(10.4,3.1){\line(0,1){1.4}}\put(11.1,1.7){\line(0,1){0.7}}
\put(11.1,3.1){\line(0,1){0.7}}
\put(11.8,1){\line(0,1){1.4}}\put(12.5,1){\line(0,1){1.4}}

\thicklines \put(11.45,2.05){\vector(1,0){0.7}}
\put(10.05,3.45){\vector(1,0){0.7}}\put(9.35,4.15){\vector(1,0){0.7}}\thinlines

\put(12.15,2.05){\line(0,1){0.7}}\put(10.75,2.05){\line(1,0){0.7}}
\put(10.75,3.45){\line(0,1){0.7}}\put(9.35,3.45){\line(1,0){0.7}}
\put(10.05,4.15){\line(0,1){0.7}}\put(8.65,4.15){\line(1,0){0.7}}

\put(12.15,2.85){\small$\theta_r$}\put(10.75,4.15){\small$\theta_3$}
\put(10.05,4.95){\small$\theta_2$}\put(10.75,1.15){\small$\theta_r$}
\put(9.65,1.85){\small$\theta_{r-1}$}\put(8.65,3.25){\small$\theta_2$}
\put(7.95,3.95){\small$\theta_1$}

\put(11.2,2.5){$\ddots$}\put(10.5,2.5){$\ddots$}

\put(1.35,4.15){\circle*{0.125}} \put(12.15,1.35){\circle*{0.125}}
\put(0.65,4.15){\line(1,0){0.7}} \put(11.45,1.35){\line(1,0){0.7}}

\put(12.85,1.0){\small $i$}\put(12.85,1.7){\small $i+1$}

\put(2.5,0){$L$ } \put(10.5,0){$\omega_i(L)$}
\end{picture}
\end{center}
\caption{$\omega_i$ acts on a Type 3 diagonal $L$}\label{D3}
\end{figure}

We next consider the case when $\lm$ has a Type 3 diagonal $L$
with content $i$. Suppose that there are $r$ boxes in $L$ and
suppose that the directions of the diagonal boxes of $L$ all go
up. Thus, the cutting strip of $\Pi$ may be written as
$[p_m,i]\uparrow[i+1,q_1]$ and the cutting strip of
$\omega_i(\Pi)$ may be written as
$[p_m,i]\blacktriangleright[i+1,q_1]$. We assume that
$\theta_{i_t}$, $1 \leq t \leq r$, is the strip in $\Pi$ that
contains the $t$-th diagonal box in $L$. Then we have
$p(\theta_{i_t})\leq i \leq q(\theta_{i_t})$, for $1\leq t \leq
r$. As illustrated in Figure \ref{D3}, under the operation of
$\omega_i$, the first strip $\theta_{i_1}=[p(\theta_{i_1}),i]$
equals
\[[p(\theta_{i_1}),q(\theta_{i_2})]=[p(\theta_{i_1}),i]\blacktriangleright[i+1,q(\theta_{i_2})].\]
If $r>1$, the last strip
\[\theta_{i_r}=[p(\theta_{i_r}),q(\theta_{i_r})]=[p(\theta_{i_r}),i]\uparrow[i+1,q(\theta_{i_r})]\]
will be cut off into $[p(\theta_{i_r}),i]$, and the other strips
\[\theta_{i_t}=[p(\theta_{i_t}),q(\theta_{i_t})]=[p(\theta_{i_t}),i]\uparrow[i+1,q(\theta_{i_t})],\]
will be twisted into
\[[p(\theta_{i_t}),q(\theta_{i_{t+1}})]=[p(\theta_{i_t}),i]\blacktriangleright[i+1,q(\theta_{i_{t+1}})],\ 2\leq t\leq r-1.\]
By Lemma \ref{inv} one may verify that $\inv(\omega_i(\Pi)) \equiv
\inv(\Pi)+r-1\ \textnormal{(mod 2)}$.

Since $i$ is the content of the ending box of $\theta_{i_1}$ and
$q(\theta_{i_t})\geq i$ for $1\leq t \leq r$, there is a unique
integer $k,\ 1\leq k\leq m-r+1$ such that $q_k=i$. So we have
$p_{k-r+1}<i+1$ and $p_{k+r}>i+1$ if $k+r<m+1$. Note in this case
$D(\Pi)$ still has the canonical form
\begin{equation}\label{H3}
(-1)^{\inv(\Pi)}\det\left[
\begin{matrix}
E&F\\
0&G
\end{matrix}\right],
\end{equation}
where
\begin{equation*}
E=\left[\begin{matrix}
s_{[p_1,\ q_1]}&\ldots&s_{[p_{k-r},\ q_1]}\\
\vdots&\ldots&\vdots\\
s_{[p_1,\ q_{k-1}]}&\ldots&
s_{[p_{k-r},\ q_{k-1}]}\\
0&\ldots&0
\end{matrix}\right]_{k\times (k-r)},
\end{equation*}

\begin{equation*}
F=\left[\begin{matrix} s_{[p_{k-r+1},\ i]\uparrow[i+1,\
q_1]}&\ldots&
s_{[p_m,\ i]\uparrow[i+1,\ q_1)]}\\
\vdots&\ldots&\vdots\\
s_{[p_{k-r+1},\ i]\uparrow[i+1,\ q_{k-1}]}&\ldots&
s_{[p_m,\ i]\uparrow[i+1,\ q_{k-1}]}\\
s_{[p_{k-r+1},\ i]}&\ldots& s_{[p_m,\ i]}
\end{matrix}\right]_{k\times (m-k+r)},
\end{equation*}

\begin{equation*}
G=\left[\begin{matrix}s_{[p_{k-r+1},\ q_{k+1}]}&\ldots&
s_{[p_m,\ q_{k+1}]}\\
\vdots&\ldots&\vdots\\
s_{[p_{k-r+1},\ q_m]}&\ldots& s_{[p_m,\ q_m]}
\end{matrix}\right]_{(m-k)\times (m-k+r)}.
\end{equation*}
On the other hand,  $D(\omega_i(\Pi))$ has the following canonical
form:
\begin{equation}\label{H4}
(-1)^{\inv(\omega_i(\Pi))}\det\left[
\begin{matrix}
E&F'\\
0&G
\end{matrix}\right]
\end{equation}
where
\begin{equation*}
F'=\left[\begin{matrix} s_{[p_{k-r+1},\
i]\blacktriangleright[i+1,\ q_1]}&\ldots&
s_{[p_m,\ i]\blacktriangleright[i+1,\ q_1)]}\\
\vdots&\ldots&\vdots\\
s_{[p_{k-r+1},\ i]\blacktriangleright[i+1,\ q_{k-1}]}&\ldots&
s_{[p_m,\ i]\blacktriangleright[i+1,\ q_{k-1}]}\\
s_{[p_{k-r+1},\ i]}&\ldots& s_{[p_m,\ i]}
\end{matrix}\right]_{k\times (m-k+r)}.
\end{equation*}
Applying Lemma \ref{llas} we get
\[s_{[p_t,\ i]\blacktriangleright[i+1,\ q_l]}=
s_{[p_t,\ i]}s_{[i+1,\ q_l]}-s_{[p_t,\ i]\uparrow[i+1,\ q_l]}\]
for $k-r+1 \leq t \leq m$ and $1 \leq l \leq k-1$. So, we may
change $F$ into $-F'$. Furthermore, by changing the signs of rows
and columns,  we obtain \eqref{H4} from \eqref{H3} with a sign
$(-1)^{r+1}$. It follows that we can transform $D(\Pi)$ into
$D(\omega_i(\Pi))$ by elementary determinantal operations. Note
that in this case there is no change on the order of the
determinants.

The other two types of diagonals can be dealt with in a similar
way. Hence we have completed the proof of Theorem \ref{thmmain}.
\qed

For example, let $\lambda=5331$, and let $\Pi_1,\Pi_2,\Pi_3,\Pi_4$
be the outside decompositions of $\lambda$ with respect to the
Jacobi-Trudi determinant, the Giambelli determinant, the
Lascoux-Pragacz determinant and the dual Jacobi-Trudi determinant.
Then we have the following relations:
$$\begin{array}{l}
\omega_{-1}\omega_{-2}\omega_{-3}(\Pi_1)=\Pi_2,\\[6pt]
\omega_1\omega_0\omega_{-1}\omega_{-2}(\Pi_2)=\Pi_3,\\[6pt]
\omega_3\omega_2\omega_{-1}\omega_{-2}(\Pi_3)=\Pi_4.
\end{array}$$
So we may transform the Jacobi-Trudi determinant for $s_{5331}$
into the Giambelli determinant for $s_{5331}$, and transform the
Giambelli determinant into the Lascoux-Pragacz determinant, and
transform the Lascoux-Pragacz determinant to the dual Jacobi-Trudi
determinant.

{\small
$$
\left|\begin{matrix}h_5&h_6&h_7&h_8\\
h_2&h_3&h_4&h_5\\
h_1&h_2&h_3&h_4\\
0&0&1&h_1
\end{matrix}\right|
\stackrel{\omega_{-3}}{\rightarrow}
\left|\begin{matrix}s_5&s_6&s_{71}\\
s_2&s_3&s_{41}\\
s_1&s_2&s_{31}
\end{matrix}\right|
\stackrel{\omega_{-2}}{\rightarrow}
\left|\begin{matrix}s_5&s_{611}&s_6\\
s_2&s_{311}&s_3\\
s_1&s_{211}&s_2
\end{matrix}\right|$$
$$
\stackrel{\omega_{-1}}{\rightarrow}
\left|\begin{matrix}s_{5111}&s_{51}&s_5\\
s_{2111}&s_{21}&s_2\\
s_{1111}&s_{11}&s_1
\end{matrix}\right|
\stackrel{\omega_{-2}}{\rightarrow}
\left|\begin{matrix}s_{51}&s_{621/1}&s_5\\
s_{21}&s_{321/1}&s_2\\
s_{11}&s_{221/1}&s_1
\end{matrix}\right|
\stackrel{\omega_{-1}}{\rightarrow}
\left|\begin{matrix}s_5&s_6&s_{71}\\
s_2&s_3&s_{41}\\
s_1&s_2&s_{31}
\end{matrix}\right|$$
$$
\stackrel{\omega_0}{\rightarrow}
\left|\begin{matrix}s_{52/1}&s_{631/2}&s_{41}\\
s_{22/1}&s_{331/2}&s_{11}\\
s_{2}&s_{31}&s_1
\end{matrix}\right|
\stackrel{\omega_1}{\rightarrow}
\left|\begin{matrix}s_{5331/22}&s_{422/11}&s_{311}\\
s_{331/2}&s_{22/1}&s_{11}\\
s_{31}&s_{2}&s_1
\end{matrix}\right|
\stackrel{\omega_{-2}}{\rightarrow}
\left|\begin{matrix}s_{422/11}&s_{42211/11}&s_{311}\\
s_{22/1}&s_{2211/1}&s_{11}\\
s_{2}&s_{211}&s_1
\end{matrix}\right|$$$$
\stackrel{\omega_{-1}}{\rightarrow}
\left|\begin{matrix}s_{311}&s_{3111}&s_{311111}\\
s_{11}&s_{111}&s_{11111}\\
s_{1}&s_{11}&s_{1111}
\end{matrix}\right|
\stackrel{\omega_2}{\rightarrow}
\left|\begin{matrix}s_2&s_{2111}&s_{21111}&s_{2111111}\\
1&s_{111}&s_{1111}&s_{111111}\\
0&s_{11}&s_{111}&s_{11111}\\
0&s_{1}&s_{11}&s_{1111}
\end{matrix}\right|
\stackrel{\omega_3}{\rightarrow}
\left|\begin{matrix}e_1&e_2&e_5&e_6&e_8\\
1&e_1&e_4&e_5&e_7\\
0&1&e_3&e_4&e_6\\
0&0&e_2&e_3&e_5\\
0&0&e_1&e_2&e_4
\end{matrix}\right|
$$
}

\section{Supersymmetric Schur Functions}

We note that Lemma \ref{llas} also holds  for the supersymmetric
Schur function $s_{\lm}=s_{\lm}(\x,\y)$ with variable set in the
form $\x=\{\ldots,x_{-1},x_0,x_1,\ldots\}$ and
$\y=\{\ldots,y_{-1},y_0,y_1,\ldots\}$. This fact can be verified
by the tableaux representation of $s_{\lm}(\x,\y)$ given in
\cite{GG} and the symmetric property of $s_{\lm}(\x,\y)$ in $\x$
and $\y$. From the Jacobi-Trudi identity for the super Schur
function $s_{\lm}(\x,\y)$ given in Goulden-Greene \cite{GG}, and
Macdonald \cite{Mac2}, one can derive a determinantal formula
analogous to Hamel and Goulden's identity for supersymmetric Schur
functions. The lattice path interpretations can also be carried
over to the supersymmetric case by using the weights given by
Chen-Yan-Yang \cite{CYY} in their study of  the super Giambelli
identity and the super Lascoux-Pragacz identity in the framework
of the  lattice path methodology. Formally speaking, we have the
following theorem.

\begin{theo}\label{thmsuper}
For any outside decomposition $\Pi$ of a skew shape $\lm$, let
$D(\Pi)$ be a determinant with the entries being the
supersymmetric ribbon Schur functions in the variable sets
$\x=\{\ldots,x_{-1},x_0,x_1,\ldots\}$ and
$\y=\{\ldots,y_{-1},y_0,y_1,\ldots\}$. Then $D(\Pi)$ equals the
supersymmetric Schur function $s_{\lm}(\x,\y)$.
\end{theo}

\vspace{.2cm} \noindent{\bf Acknowledgments.} We would like to
thank Professor Alain Lascoux for valuable comments.  This work
was done under the auspices of the 973 Project on Mathematical
Mechanization of the Ministry of Science and Technology, and the
National Science Foundation of China.

\end{document}